\documentstyle[12pt]{article}
\begin{document}
\title{A Note on $n$-ary Poisson Brackets}
\author{{\normalsize by}\\Peter W. Michor and Izu Vaisman}
\date{}
\maketitle
{\def\thefootnote{*}\footnotetext[1]%
{{\it 1991 Mathematics Subject Classification} 58 F 05. 
\newline\indent PWM was supported by `Fonds zur F\"orderung der
wissenschaftlichen Forschung, Projekt P~10037~PHY'.
During the work on this paper, IV was a visitor at the Erwin Schr\"odinger 
International Institute for Mathematical Physics
in Vienna, Austria, and he expresses here his thanks to ESI for 
invitation and support.}}
\begin{center} \begin{minipage}{12cm}
A{\footnotesize BSTRACT.
A class of $n$-ary Poisson structures of constant rank is indicated.
Then, one proves that the ternary 
Poisson brackets are exactly those which are defined by
a decomposable $3$-vector field. The key point is the 
proof of a lemma 
which tells that an $n$-vector $(n\geq3)$ is decomposable
iff all its contractions with up to $n-2$ covectors are decomposable.}
\vspace{5mm}
\end{minipage} \end{center}

\indent 
In the last years, several authors have studied 
generalizations of Lie  algebras to various types of n-ary algebras, e.g.,
\cite{{Fp},{Tk},{HW},{MV},{VV}}.  In the same time, and
intended  to physical applications, the new types of  algebraic structures 
were considered in the case of the algebra $C^{\infty}(M)$ of  
functions on a $C^{\infty}$ manifold $M$, under the assumption that 
the operation is a  
derivation of each entry separately. In this way one got the Nambu-Poisson 
brackets, e.g., \cite{{Tk},{Gt},{AG},{F}}, and the generalized Poisson 
brackets \cite{{Az1},{Az2}}, etc. In this note, we 
write down the
characteristic conditions of the $n$-ary
generalized Poisson structures in a new form,
and give an example of an $n$-ary 
structure of constant rank $2n$, for any $n$ even or odd. Then, we prove
that the ternary Poisson brackets
are exactly the brackets defined by the decomposable $3$-vector fields.
The key point in the proof of this result is a 
lemma (that seems to appear also in \cite{W}),
which tells that an $n$-vector $P$
is decomposable iff $i(\alpha_{1})...i(\alpha_{k})P$ is decomposable,
for any choice of covectors $\alpha_{1},...,\alpha_{k}$, 
where $k$ is fixed, and such that $1\leq k
\leq n-2$.
 
Our framework is the $C^{\infty}$ category. If $M$ is an 
$m$-dimensional manifold, an {\em $n$-ary Poisson bracket or 
structure} (called generalized Poisson structure in \cite{{Az1},{Az2}}),
with the {\em Poisson $n$-vector or tensor} $P$, is a bracket of the form 
$$\{f_{1},...,f_{n}\}=P(df_{1},...,df_{n}) \hspace{5mm} 
(f_{1},...,f_{n}\in C^{\infty}(M)), \eqno{(1)}$$ where 
$P\in\Gamma\wedge^{n}TM$ is an $n$-vector (i.e., a completely 
skew-symmetric contravariant tensor) field, and the following 
{\em generalized Jacobi identity} of order $n$ \cite{MV} is satisfied 
$$\sum_{\sigma\in S_{2n-1}}(sign\,\sigma)
\{\{f_{\sigma_{1}},...,f_{\sigma_{n}}\}, 
f_{\sigma_{n+1}},...,f_{\sigma_{2n-1}}\}=0,\eqno{(2)}$$
$S_{2n-1}$ being the symmetric group.

For $n=2$ the bracket is a usual Poisson bracket e.g., \cite{V2}.
In this note, we always assume $n\geq3$.
\proclaim 1 Proposition. 
The $n$-vector field $P\in\Gamma\wedge^{n}TM$  defines an $n$-ary
Poisson bracket
iff either $n$ is even and the Schouten-Nijenhuis bracket 
$[P,P]=0$, or $n$ is odd and $P$ satisfies the conditions 
$$(i(\alpha)P)\wedge(i(\beta)P)=0\hspace{5mm}\forall\alpha,\beta\in 
T^{*}M,\eqno{(A)}$$
$$\sum_{u=1}^{m}(i(dx^{u})P)\wedge
(L_{\partial/\partial x^{u}}P)=0,\eqno{(D)}$$
where $(x^{u})$ are local coordinates on $M$, and $L$ denotes Lie 
derivative.  \par 
\noindent{\bf Proof.} The left hand side of (2) contains only first 
and second order 
derivatives, and is skew symmetric in the arguments $f_{i}$. Hence,
to ensure (2) it is enough to 
ask it to hold for the case of the local functions $f_{i}=x^{a_{i}}$, 
$(i=1,...,2n-1)$, and for the case of the functions 
$f_{1}=x^{u}x^{v}$,  $f_{i}=x^{a_{i}}$, $(i=2,...,2n-1)$, at $x=0$. In the 
first case the result is
$$\sum_{u=1}^{n}P^{u[a_{1}...a_{n-1}}\frac{\partial}{\partial x^{u}}
P^{a_{n}...a_{2n-1}]}=0,\eqno{(3)}$$
and in the second case the result is
$$P^{v[a_{2}...a_{n}}P^{a_{n+1}...a_{2n-1}]w}+ 
P^{w[a_{2}...a_{n}}P^{a_{n+1}...a_{2n-1}]v}=0,\eqno{(4)}$$
where square brackets denote index alternation. 
Now, (4) is equivalent to (A) if $n$ is odd, and it is an identity if $n$ is
even. Then, (3) is equivalent to (D), and the use of the coordinate 
expression of the 
Schouten-Nijenhuis bracket (e.g., \cite{V2}) shows that, for $n$ even, (D)
is equivalent to $[P,P]=0$. Q.e.d.
 
We call (A) and (D) the {\em algebraic} and the {\em differential
condition}, respectively.
The coordinate expressions (3), (4), and their equivalence with 
$[P,P]=0$ in the $n$-even case, were also established in 
\cite{{Az1},{Az2}}. In the $n$-odd case, the differential condition (D)
has no independent invariant meaning, 
and it must be associated with the algebraic condition (A).  

It is also important to notice that, since (A,D) always hold at 
the zeroes of $P$, $P$ defines an $n$-ary Poisson bracket iff it does 
so on the  subset $U\subseteq M$ where $P\neq0$.\\

Before going on, we need some general facts about $n$-vectors 
$P\in\wedge^{n}L$, 
where $L$ is an $m$-dimensional (e.g., real) linear space. 
First, $P$ defines a linear mapping
$\sharp_{P}:\wedge^{n-1}L^{*}\rightarrow L$ given by
$$\sharp_{P}(\lambda)=i(\lambda)P,\hspace{5mm}\lambda\in\wedge^{n-1}L^{*}.
\eqno{(5)}$$
We will say that $rank\,\sharp_{P}=dim\,im\,\sharp _{P}$ is the 
{\em rank} of $P$.  
(This definition is equivalent with the one used in older books on
exterior algebra e.g., \cite{Sz}, which refered to $L^{*}$ rather than $L$,
and where the vectors of $im\,\sharp_{P}$ were seen as the right hand 
side of the equations of the {\em adjoint system} of $P$.)

If $rank\,P=dim\,L$, we say that $P$ is {\em non degenerate} 
({\em regular} in \cite{Sz}).
On the other hand, 
if $rank\,P=n$, $P$ is {\em decomposable} i.e., 
there are vectors 
$W_{a}\in L$ $(a=1,...,n)$ such that $P=W_{1}\wedge...\wedge W_{n}$.
We recall the existence of classical decomposability conditions known
as the {\em Pl\"ucker conditions} e.g. \cite{GH},
which we will write down later, in  the proof of Lemma 3.

The skew symmetry of $P$ implies $im\,\sharp_{P}=Ann\,(A(P))$,
where $A(P):=\{\alpha\in L^{*}\,/\,i(\alpha)P=0\}$, hence, 
$rank\,P=m-dim\,A(P)$. Notice also that 
$P\in\wedge^{n}(im\,\sharp_{P})$. Indeed, if 
$L=im\,\sharp_{P}\oplus K$, we have an expression 
$$P=\sum_{u+v=n}Q_{u}\wedge S_{v},
\hspace{5mm}Q_{u}\in\wedge^{u}(im\,\sharp_{P}),\,
S_{v}\in\wedge^{v}K,$$
and, since $K^{*}$ may be identified with $A(P)$, 
the previous expression reduces to $P=Q_{n}$. On the other hand,
if $P\in\wedge^{n}U$, where
$U$ is a subspace of $ L$, $im\,\sharp_{P}\subseteq U$.
Therefore, $im\,\sharp_{P}$
is the minimal 
subspace $S$ of $L$ such that $P\in\wedge^{n}S$.

The $n$-vector $P$ will be called {\em irreducible} if there is no 
decomposition $im\,\sharp_{P}=S_{1}\oplus S_{2}$ 
where $dim\,S_{1}=n$, and where $P=P_{1}+P_{2}$
with $0\neq P_{1}\in\wedge^{n}S_{1}$, $0\neq P_{2}\in\wedge^{n}S_{2}$.
If such a decomposition exists, $P$ is {\em reducible}, and, because
$\forall\lambda\in\wedge^{n-1}L^{*}$, 
$i(\lambda)P=i(\lambda)P_{1}+i(\lambda)P_{2}$,
we have $S_{1}=im\,\sharp_{P_{1}}$, $S_{2}=im\,\sharp_{P_{2}}$.
 From these definitions, it follows that any $n$-vector $P$ may 
be (non uniquely) written
under the form
$$P=\sum_{i=0}^{s-1}V_{in+1}\wedge...\wedge V_{in+n}+P',\eqno{(6)}$$
where $V_{a}$ $(a=1,...,sn)$ are independent vectors and 
$P'\in\wedge^{n}U$, where $U$ is a complement of $span\,\{V_{a}\}$ in $L$,
is irreducible with
$rank\,P'=rank\,P-sn$.

An $n$-vector $P$, which satisfies
condition (A) $\forall\alpha,\beta\in L^{*}$,
must be irreducible since otherwise, and with the notation above,
we have $(i(\alpha)P)\wedge(i(\beta)P)\neq0$ if $i(\alpha)P_{1}\neq0$,
$i(\alpha)P_{2}=0$, $i(\beta)P_{1}=0$, $i(\beta)P_{2}\neq0$. Of course, if
$rank\,P<2n-2$,
$P$ is irreducible and (A) is an identity.

Now, we come back to the manifold $M$. Then, if $P\in\Gamma\wedge^{n}TM$,
$rank\,P$ is a lower semicontinuous function on $M$.
The following Proposition gives an interesting class of Poisson, but
not Nambu-Poisson, generally,
$n$-vectors, for an arbitrary even or odd order $n$.
A first example of an $n$-ary Poisson
structure of an even order $n$ was given in \cite{Az1}, and it was a linear 
structure on the dual of a simple Lie algebra. We know of no previous
examples of $n$-ary Poisson structures which are not Nambu-Poisson
structures. 

\proclaim 2 Proposition. Assume that $P\in\Gamma\wedge^{n}TM$,
and that $\forall x\in M$ there is an open neighbourhood $U_{x}$
such that $P/_{U_{x}}$ can be written as
$$P=\frac{1}{h!(n-h)!}\sum_{\sigma\in S^{n}}(sign\,\sigma)
V_{\sigma_{1}}\wedge ... \wedge V_{\sigma_{h}}
\wedge W_{\sigma_{h+1}}\wedge...\wedge W_{\sigma_{n}},
\eqno{(7)}$$
where $S_{n}$ is the permutation group, 
$(V_{i},W_{j})$ $(i,j=1,...,n)$ are independent vector fields on 
$U_{x}$, and $h$ is a fixed integer such that $0\leq2h\leq n-3$.
Then $P$ is a Poisson $n$-vector of constant rank, equal to $2n$
if $h\neq0$, and to $n$ if $h=0$.  

\noindent{\bf Proof.} The condition on $h$ was chosen such that the 
left hand sides of (A) and (D), as well as $[P,P]$, consist of sums 
of wedge products where at least one of the vectors $W_{j}$ must be 
wedge multiplied by itself. If $h=0$, $P=W_{1}\wedge...\wedge W_{n}$ 
i.e., $P$ is decomposable and of rank $n$.
If $h\neq0$, by using a basis of vectors that starts with 
$(V_{i},W_{j})$ and its dual cobasis in order to
obtain generators of $im\,\sharp_{P}$, we see that 
$im\,\sharp_{P}=span\{V_{i},W_{j}\}$ hence, $rank\,P=2n$.
(The particular case of an even order, decomposable $n$-vector $P$
was noticed in \cite{Az2}.)
Q.e.d.

We will say that a tensor
$P$ of form (7) with $h\neq0$
is a {\em semi-decomposable $n$-vector}.
An $n$-ary Poisson structure (bracket) defined by a (semi-) decomposable
$n$-vector field will be called a {\em (semi-)decomposable $n$-ary Poisson
structure (bracket)}. 

\proclaim 3 Lemma. Let $L$ be an $m$-dimensional vector space, and 
$P\in\wedge^{n}L$. 
Then $P$ is decomposable iff for any fixed number $k$ such that
$1\leq k\leq n-2$, and any set of $k$ covectors
$\alpha_{u}\in L^{*}$ $(u=1,...,k)$, the $(n-k)$-vector
$i(\alpha_{1})...i(\alpha_{k})P$ is decomposable.\par

Appearently,
this lemma is included in formula (4), page 116 
of \cite{W}. Our proof is different. 

\noindent{\bf Proof.} 
It is enough to prove the result for $k=1$, 
and, on the other hand, the proof and the
result do not hold for $k>n-2$.

As already recalled,
decomposability of $P$ is characterized by the Pl\"ucker conditions.
These may be written in one of the following equivalent forms 
\cite{GH}:
$i)$ $\forall V\in im\,\sharp_{P},$ $V\wedge P=0$, $ii)$ 
$\forall\lambda\in\wedge^{n-1}L^{*},$ $(i(\lambda)P)\wedge P=0$.
(The first form of the conditions is rather obvious, and the second is 
equivalent since $i(\lambda)P$ are exactly the vectors of 
$im\,\sharp_{P}$.) 

 From 
$$i(\alpha)[(i(\mu)i(\alpha)P)\wedge P]
=-[i(\mu)(i(\alpha)P)]\wedge(i(\alpha)P
\hspace{5mm}(\alpha\in L^{*},\mu\in\wedge^{n-2}L^{*}),$$
we see that if $P$ satisfies condition $ii)$, i.e.,  if $P$ is 
decomposable, so are all $i(\alpha)P$, $\alpha\in L^{*}$.

Now, assume that 
all $i(\alpha)P$ are decomposable, and take $\epsilon^{1}\in L^{*}$
such that $i(\epsilon^{1})P\neq0$. Then, $\exists e_{a}\in ker\,
\epsilon^{1}\subseteq L$
$(a=2,...,n)$ such that
$$i(\epsilon^{1})P=e_{2}\wedge...\wedge e_{n}.\eqno{(8)}$$
Let us also take $e_{1}\in L$ such that $\epsilon^{1}(e_{1})=1$, and 
denote by $L_{1}$ the $n$-dimensional subspace
$span\{e_{1},...,e_{n}\}$ of $L$, and by
$L_{2}$ an arbitrary complement of $span\{e_{2},...,e_{n}\}$ in 
$ker\,\epsilon^{1}$. 
Then $L=L_{1}\oplus L_{2}$, and we have an expression
$$P=\rho e_{1}\wedge...\wedge e_{n}+\sum_{i=1}^{n-1}P'_{i}\wedge
P''_{i}+P''_{n},\eqno{(9)}$$
where $\rho\in{\bf R}$, $P'_{i}\in\wedge^{n-i}L_{1}$,
$P''_{i}\in\wedge^{i}L_{2}$,
$P''_{n}\in\wedge^{n}L_{2}$. Moreover, (8) implies $\rho=1$ and
$$i(\epsilon^{1})P'_{i}=0\hspace{5mm}(i=1,...,n-1).\eqno{(10)}$$
(If some $P''_{i}=0$ we will also assume $P'_{i}=0$.)

Let $\epsilon^{i}\in L^{*}$ be covectors which vanish on $L_{2}$, and 
are such that $\epsilon^{i}(e_{j})=\delta_{j}^{i}$ $(i,j=1,...,n)$. 
According to our hypothesis,
the $(n-1)$-vectors 
$$i(\epsilon^{a})P=(-1)^{a-1}e_{1}\wedge...
\wedge\hat e_{a}\wedge...\wedge e_{n}
+\sum_{i=1}^{n-1}(i(\epsilon^{a})P'_{i})\wedge P''_{i}$$ 
$(a=2,...,n)$, where the hat denotes the absence of the factor, must 
also be decomposable. In view of (10),
for $$\lambda=\epsilon^{1}\wedge...\wedge\hat \epsilon^{a}\wedge...\wedge
\hat\epsilon^{b}\wedge...\wedge\epsilon^{n}\hspace{5mm}(b\neq a),$$ we have
$i(\lambda)i(\epsilon^{a})P=\pm e_{b}$, where $b=2,...,n$, and the sign
depends
on whether $a<b$ or $b<a$,
and the
Pl\"ucker condition $ii)$ yields
$$e_{b}\wedge (i(\epsilon^{a})P)=
\sum_{i=1}^{n-1}e_{b}\wedge (i(\epsilon^{a})P'_{i})\wedge
P''_{i}=0.\eqno{(11)}$$
This implies $e_{b}\wedge(i(\epsilon^{a})P'_{i})=0$, $i=1,...,n-1$,
and the $(n-i-1)$-vector $i(\epsilon^{a})P'_{i}$ 
belongs to the ideal generated by $e_{2}\wedge...\wedge \hat
e_{a}\wedge...\wedge e_{n}$.
Therefore,
$i(\epsilon^{a})P'_{i}=0$, except for $i=1$, and, using again (10),
$$i(\epsilon^{a})P'_{1}=\kappa e_{2}\wedge...\wedge\hat e_{a}\wedge
...\wedge e_{n}\hspace{5mm}(\kappa\in{\bf R}).$$
Accordingly,
$$P'_{1}=(-1)^{a-1}\kappa e_{2}\wedge...\wedge e_{a}\wedge...\wedge
e_{n},\,
P'_{2}=0,\,...,\,P'_{n-1}=0,$$
and we deduce
$$P=e_{2}\wedge...\wedge
e_{n}\wedge((-1)^{n-1}e_{1}+(-1)^{a-1}P''_{1})+P''_{n}.
\eqno{(12)}$$
In other words, $P$ is reducible. But, then, if we take $\alpha=
\beta+\gamma\in L^{*}$, where $\beta$ vanishes on the second term of 
(12) but not on the first,
and $\gamma$ vanishes on the first term but not on the second, we see that
$i(\alpha)P$ is not decomposable unless $P''_{n}=0$. Hence, our
$P$ must be decomposable. Q.e.d.

The decomposable $n$-vectors $P\neq0$ are important because they define
the $n$-planes, via $im\,\sharp_{P}$. $\forall\alpha\in L^{*}$, one has
$$im\,\sharp_{i(\alpha)P}\subseteq(ker\,\alpha)
\cap(im\,\sharp_{P}),\eqno{(13)}$$
while, if $i(\alpha)P\neq0$, the subspaces in the right hand side of 
(13) are transversal in $L$, and the intersection has the dimension 
$rank\,P-1$. Hence, if $rank\,(i(\alpha)P)=rank \,P-1$ one has
$$im\,\sharp_{i(\alpha)P}=(ker\,\alpha)\cap(im\,\sharp_{P}),\eqno{(14)}$$
In particular, this is always true if $P$ is decomposable.

Notice that, because of (6),  if $P$ is reducible, $\exists\alpha\in L^{*}$
such that $i(\alpha)P\neq0$ and 
$rank\,(i(\alpha)P)\leq rank \,P-n<rank \,P-1$. 
Therefore, if (14) holds $\forall\alpha\in L^{*}$ with $i(\alpha)P\neq0$
hence, $rank\,(i(\alpha)P)=rank \,P-1$, $P$ is irreducible.

Now, as a particular case of Lemma 3 we get
\proclaim 4 Lemma. Let $L$ be an $m$-dimensional linear space,
and $P\in\wedge^{3}L$. Then, if 
$P$ satisfies condition (A), $\forall\alpha,\beta\in L^{*}$, 
$P$ is decomposable. \par 
\noindent{\bf Proof.} For any 
$\alpha\in L^{*}$, the bivector $Q=i(\alpha)P$ is decomposable,
since condition (A) implies $Q\wedge Q=0$, and in the case of a bivector
this is equivalent with the
Pl\"ucker decomposability condition $ii)$ above. Indeed,
if $Q$ is decomposable, obviously $Q\wedge Q=0$. 
Conversely, $Q\wedge Q=0$ implies 
$$i(\alpha)(Q\wedge Q)=2i(\alpha)Q\wedge Q=0,
     \hspace{5mm}\forall\alpha\in L^{*}.$$ 
Therefore, the result follows from Lemma
3. Q.e.d.

 From Proposition 2 and Lemma 4 we get
\proclaim 5 Theorem. A $3$-vector field $P$ defines a ternary Poisson 
bracket on the manifold $M$ iff, around every point $x\in M$, $P$ is 
decomposable. 
\par
In particular, there is no differential condition to be imposed on
a Poisson trivector, since condition (D) is a consequence
of decomposability.\par
Let us also notice 
\proclaim 6 Corollary. A ternary Poisson bracket is a Nambu-Poisson
bracket iff the distribution $im\,\sharp_{P}$ is involutive.\par
This follows from the well known fact that, where non zero, a 
Nambu-Poisson bracket is a Jacobian determinant e.g., \cite{{AG},{Gt}}.

We finish by a few more related remarks. 

First, it is known that all the Nambu-Poisson tensors are decomposable.
This follows from the fact that they must satisfy the algebraic condition
\cite{Tk}
$$\sum_{k=1}^{n}[P^{b_{1}\ldots b_{k-1}ub_{k+1}\ldots b_{n}} 
P^{va_{2}\ldots a_{n-1}b_{k}}+
P^{b_{1}\ldots b_{k-1}vb_{k+1}\ldots b_{n}}
P^{ua_{2}\ldots a_{n-1}b_{k}}]=0.\eqno{(N1)}$$
In \cite{AG} there is an algebraic proof of the fact that (N1)
implies decomposability. Lemma 3 above allows for a very short
proof of the same result. Namely, (N1) is equivalent to
$$
i(\alpha)P\wedge i(\Phi)i(\beta)P + i(\beta)P\wedge i(\Phi)i(\alpha)P 
= 0, \eqno{(N2)}$$
$\forall \alpha, \beta \in T^*M$ and
$\forall \Phi\in \wedge^{n-2} T^*M 
$,
and (N2) is the polarization of the, once more equivalent, condition
$$
i(\alpha)P\wedge i(\Phi)i(\alpha)P 
= 0\hspace{5mm} \forall \alpha\in T^*M, 
\forall \Phi\in \wedge^{n-2} T^*M.\eqno{(N3)}
$$
By the Pl\"ucker relations this means that $i(\alpha)P$ is 
decomposable for all $\alpha\in T^*M$, which by Lemma
3 is 
equivalent to $P$ being decomposable. 

The above proof clarifies the relation between the Pl\"ucker
and the Nambu decomposability conditions.

Second, an $n$-ary Poisson structure $P$ of constant rank defines a
tensorial 
$G$-structure on the manifold $M$, and it is natural to ask what are 
the integrability conditions of this structure. 
Following are two examples of integrable
structures, written by means of the corresponding systems of 
local coordinates:
$$P=\frac{1}{h!(n-h)!}\sum_{\sigma\in S_{n}}(sign\,\sigma)
\frac{\partial}{\partial y^{\sigma_{1}}}
\wedge ... \wedge \frac{\partial}{\partial y^{\sigma_{h}}}
\wedge \frac{\partial}{\partial x^{\sigma_{h+1}}}
\wedge...\wedge
\frac{\partial}{\partial x^{\sigma_{n}}}, \eqno{(15)}$$
$$P=\sum_{i=0}^{s-1}\frac{\partial}{\partial x^{2iu+1}}\wedge...\wedge 
\frac{\partial}{\partial x^{2(iu+u)}}.\eqno{(16)}$$
In (16), we have an integrable Poisson $n$-vector where $n=2u$ is 
even, and the $n$-vector 
is reducible and does not satisfy condition (A). We may say that $P$ 
of (16) is the generalization of a symplectic structure since the 
latter can be defined by the same formula for $n=2$.
 
The last remark is that, if a tensor field $P$ defines a
(semi-) decomposable
$n$-ary Poisson structure on a manifold $M$, there is an 
interesting Grassmann subalgebra on $M$, namely, 
$\Sigma M=\Sigma^{q}M$ where 
$$\Sigma^{q}M:=\{U\in\Gamma\wedge^{q}TM\;/\;\forall\alpha\in T^{*}M,\; 
(i(\alpha)P)\wedge(i(\alpha)U)=0\}.\eqno{(17)}$$
On $\Sigma M$, a differential operator 
$\delta:\Sigma^{q}M\rightarrow\Sigma^{q+n-1}M$ may be defined by
$$\delta U=\sum_{u=1}^{m}((i(dx^{u})P)\wedge(L_{\partial/\partial x^{u}}U)+ 
(i(dx^{u})U)\wedge(L_{\partial/\partial x^{u}}P)),\eqno{(18)}$$
where $(x^{u})$ are local coordinates on $M$.
The fact that $\delta$ is well
defined follows by the same argument as in the proof of Proposition 2,
while using the polarization of condition
(17) with respect to $\alpha$.
 
In particular, $P\in\Sigma^{n}M$ and $\delta P=0$, 
and for $\forall f\in C^{\infty}(M)$,
$f\in\Sigma^{0}M$ and $\delta f=i(df)P$. 
Generally, we do not have $\delta^{2}=0$, and only 
the {\em  twisted cohomology} 
$ker\delta/(im\delta\cap ker\delta)$ \cite{V}
can be considered. 
{\small P. W. Michor, 
Institut f\"ur Mathematik, Universit\"at Wien,
Strudlhofgasse 4, A-1090 Wien, Austria, {\it and}:
Erwin Schr\"odinger Institute, Boltzmanngasse 9, A-1090, Wien,
Austria. E-mail: michor@esi.ac.at}\\
{\small I. Vaisman: Dept. of Math., Univ. of Haifa,
Haifa 31905, Israel.} \\{\small E-mail: vaisman@math.haifa.ac.il}\\

\begin{thebibliography}{xx} 
\bibitem{AG} D. Alekseevsky and P. Guha, On decomposability of 
Nambu-Poisson Tensors. Acta Math. Univ. Commenianae, 65 (1996), 1-10.
\bibitem{Az1} J. A. de Azc\'arraga, A. M. Perelomov and J. C. P\'erez Bueno,
New generalized Poisson structures. J. Phys. A: Math. Gen. 29 (1996),
L151-157.
\bibitem{Az2} J. A. de Azc\'arraga, J. M. Izquierdo and J. C. P\'erez Bueno,
On the higher order generalization of Poisson structures. Preprint
hep-th/97033019 v2. J. Phys. A, to appear.
\bibitem{F} G. Dito, M. Flato, D. Sternheimer and L. Takhtajan, Deformation
quantization and Nambu Mechanics. Comm. in Math. Phys., 183 (1997), 1-22.
\bibitem{Fp} V. T. Filippov, $n$-ary Lie algebras. Russian Sibirskii 
Math. J., 24 (1985), 126-140.
\bibitem{Gt} Ph. Gautheron, Some remarks concerning Nambu mechanics. 
Lett. in Math. Phys., 37 (1996), 103-116.
\bibitem{GH} Ph. Griffiths and J. Harris, Principles of Algebraic 
Geometry, 2nd. ed. J. Willey \& Sons, New York, 1994.
\bibitem{HW} P. Hanlon and M. Wachs, On Lie $k$-algebras. Adv. in 
Math., 113 (1995), 206-236.
\bibitem{Sz} W. \'Slebodzi\'nski, Exterior forms and their applications. 
PWN-Polish Scientific Publishers, Warszawa, 1970.
\bibitem{MV} P. W. Michor and A. M. Vinogradov, $n$-ary Lie and associative
algebras. Rend. Sem. Mat. Torino, 53 (1996), 373-392.
\bibitem{Tk} L. Takhtajan, On Foundations of Generalized Nambu Mechanics.
Comm. Math. Phys., 160 (1994), 295-315.
\bibitem{V} I. Vaisman, New examples of twisted cohomologies, Bollettino
Unione Mat. Italiana, 7-B (1993), p. 355-368.
\bibitem{V2} I. Vaisman, Lectures on the Geometry of Poisson 
Manifolds. Progress in Math., vol. 118, Birkh\"auser Verlag, Basel 1994.
\bibitem{VV} A. Vinogradov and M. Vinogradov, On multiple 
generalizations of Lie algebras and Poisson manifolds. Contemporary 
Math., 219 (1998), 273-287. 
\bibitem{W} R. Weitzenb\"ock, Invariantentheorie. P. Noordhoff, 
Groningen, 1923.
\end{thebibliography}
\end{document}